\documentclass[11pt,leqno]{article}

\usepackage{amsmath}
\usepackage{amssymb}
\usepackage{amscd}
\usepackage{graphicx}
\input{epsf}

\newtheorem{MainTheorem}{Theorem}

\newtheorem{Theorem}{Theorem}[section]

\newtheorem{Lemma}[Theorem]{Lemma}
\newtheorem{Proposition}[Theorem]{Proposition}

\newtheorem{Example}[Theorem]{Example}

\,
\,
\,

\def\qed{\hfill $\bullet$}
\def\proof{\noindent {\em Proof:}\ }

\def\set-up{\noindent {\em Set-Up:}\ }

\def\ra{\rightarrow}

\def\<{\langle}
\def\>{\rangle}

\def\"{\,^{\prime\prime}}

\def\iff{\Leftrightarrow}

\begin{document}

\title{On Shift Dynamics for Cyclically Presented Groups}

\author{William A. Bogley}
\date{March 6, 2013}

\maketitle

\begin{abstract}
For group presentations with cyclic symmetry, there is a connection between asphericity and the dynamics of the shift automorphism. For the class of groups $G_n(k,l)$ described by the cyclic presentations $\mathcal{P}_n(k,l) = (x_i:x_ix_{i+k}x_{i+l}\ (i \mod n))$ and studied extensively by G.~Williams and M.~Edjvet \cite{EdjvetWilliams}, the shift acts freely on the nonidentity elements of $G_n(k,l)$ if and only if the presentation $\mathcal{P}_n(k,l)$ is combinatorially aspherical in the sense of \cite{CCH}. The shift has a nonidentity fixed point precisely when $G_n(k,l)$ is finite. \textbf{MSC} (2010) 20F05, 20E36. \textbf{Keywords:} cyclically presented group, shift, dynamics, asphericity.
\end{abstract}

\section{Cyclically Presented Groups}

Since the introduction of the Fibonacci groups by J.~C.~Conway in 1965 \cite{JC,JC+}, group presentations that admit cyclic symmetry have provoked innovative approaches in combinatorial group theory \cite{JWW,N,TRevisit} and have manifested interesting connections with three-manifold geometry and topology \cite{CHK,CRS,HKM,MV,Telloni}, number theory \cite{Mac,O}, and computational group theory \cite{EdjvetIrred,H,HRS}. By their very nature, these groups come pre-equipped with a periodic automorphism, called the shift, that cyclically permutes the generators. The purpose of this paper is to describe some connections between the dynamics of the shift and the algebraic, combinatorial, and topological aspects of these groups.

%Let $X$ be a basis for a free group $F = F(X)$ and let $\Theta$ be a group of permutations of $X$. We can then view $\Theta$ as a group of automorphisms of $F$. Given a word $w = w(X)$ representing an element of $F$, we then have the \textbf{symmetric group presentation} $\mathcal{P}(X,\Theta,w)$ with generating set $X$ and relator set consisting of the orbit of $w$ under the action of $\Theta$ on $F$. The permutation group $\Theta$ now acts by automorphisms on the group presented by $\mathcal{P}(X,\Theta,w)$. We focus on the case of a cyclic permutation group generated by a transitive permutation of a finite set $X$.

Given a positive integer $n$, the \textbf{shift} on the free group with basis $\{x_0,\ldots,x_{n-1}\}$ is the length-preserving automorphism $\theta$ given by $\theta(x_i) = x_{i+1}$ with subscripts modulo $n$. Given a cyclically reduced word $w = w(x_0,\ldots,x_{n-1})$, there is the \textbf{cyclic group presentation}

$$
\mathcal{P}_n(w) = \left(x_0,\ldots,x_{n-1}\ : w, \theta(w), \ldots, \theta^{n-1}(w)\right)
$$
that defines the \textbf{cyclically presented} group $G_n(w)$. The shift can be viewed as an automorphism $\theta \in \mathrm{Aut}(G_n(w))$ having exponent $n$. D.~L.~Johnson \cite{JExt} and others have shown that the shift has order $n$ in many interesting cases, which amounts to saying that the cyclic group $C_n$ of order $n$ acts faithfully on $G_n(w)$ via the shift.

Our main results deal with fixed points of the shift and its powers. A cyclic presentation $\mathcal{P}_n(w)$ is \textbf{orientable} if $w$ is not a cyclic permutation of the inverse of any of its shifts. Lemma \ref{orientableShape} below states that $\mathcal{P}_n(w)$ is orientable unless $n = 2m$ is even and there is a word $u = u(x_0,\ldots,x_{n-1})$ such that $w = u\theta^m(u)^{-1}$. Note that $u$ is fixed by $\theta^m \in \mathrm{Aut}(G_{2m}(u\theta^m(u)^{-1}))$.

\vspace{.25in}

%\noindent \textbf{Theorem A}\ \ \textit{If the cyclic presentation $\mathcal{P}_n(w)$ is orientable and combinatorially aspherical, then the cyclic group of order $n$ acts freely via the shift on the nonidentity elements of $G_n(w)$.}

\noindent \textbf{Theorem A}\ \ \textit{If $\mathcal{P}_n(w)$ is orientable and combinatorially aspherical, then $C_n$ acts freely via the shift on the nonidentity elements of $G_n(w)$.}

\vspace{.25in}

\noindent Proposition \ref{CA} contains a working definition of combinatorial asphericity in terms of identity sequences for group presentations \cite{CCH}. Theorem A says that when the cyclic presentation $\mathcal{P}_n(w)$ is orientable and combinatorially aspherical, the orbits of the shift partition the nonidentity elements of $G_n(w)$ into pairwise disjoint $n$-cycles. In particular, $\theta$ and all of its powers are fixed-point free and $\theta$ has order $n$. Another way to view this is that if some power of the shift has a nonidentity fixed point, then the presentation $\mathcal{P}_n(w)$ must support an interesting spherical diagram \cite{CCH}. Perhaps one could prove Theorem A in this way, but that is not the path followed here.

%We rely on the well-known fact that the shift on $G_n(w)$ is realized by an inner automorphism of the semidirect product $E_n(w) = G_n(w) \rtimes_\theta C_n$, which admits a two-generator two-relator presentation of the form $(a,x:a^n,W(a,x))$. This in turn determines a relative presentation $(C_n,x:W)$ where $W$ is obtained by interpreting $W(a,x)$ in the free product $C_n \ast \<x\>$. Now $G_n(w)$ is the kernel of a retraction of $E_n(w)$ onto $C_n$, but depending on $n$ and $w$, there can be several different retractions $\nu^f:E_n(w) \ra C_n$, each determined by the image $\nu^f(x) = a^f \in C_n$. We use a rewriting process $\rho^f$ to describe the kernel of any such retraction as the group given by the cyclic presentation $\mathcal{P}_n(\rho^f(W))$ (Theorem \ref{splitCrit}). The word $w$ is recovered as $\rho^0(W) = w$. As $f$ ranges through allowable values, which are defined modulo $n$, this leads to commensurable families of cyclically presented groups. These commensurable groups need not be isomorphic, but they do have identical dynamics under their respective shifts (Lemma \ref{TransferDynamic}). Commensurability status can also be an aid in identifying the structure of $G_n(w)$.

Many have noted that the shift on $G_n(w)$ is realized by an inner automorphism of the semidirect product $E_n(w) = G_n(w) \rtimes_\theta C_n$, which admits a two-generator two-relator presentation of the form $(a,x:a^n,W(a,x))$. This in turn determines a relative presentation $(C_n,x:W)$ where $W$ is obtained by interpreting $W(a,x)$ in the free product $C_n \ast \<x\>$. Now $G_n(w)$ is the kernel of a retraction of $E_n(w)$ onto $C_n$, but depending on $n$ and $w$, there can be several different retractions $\nu^f:E_n(w) \ra C_n$, each determined by the image $\nu^f(x) = a^f \in C_n$. A rewriting process $\rho^f(W)$ is introduced in Section \ref{sec:Retrations} to show that the kernel of such a retraction is cyclically presented (Theorem \ref{splitCrit}). The original word $w$ is recovered as (a shift of) $\rho^0(W)$. As $f$ ranges through allowable values there arise families of commensurable cyclically presented groups. These commensurable groups need not be isomorphic, but they do have identical dynamics under their respective shifts (Lemma \ref{TransferDynamic}). Commensurability status can also be an aid in identifying the structure of $G_n(w)$ (Lemma \ref{lemma:metacyclic}).

The theory of aspherical orientable relative presentations \cite{BP} aids in the proof of Theorem A and is the focus of its practical applications. The cellular model of a relative presentation is described in Section \ref{sec:Asphericity}. (See also \cite[Section 4]{BP}.)

\vspace{.25in}

\noindent \textbf{Theorem \ref{freeAction}}\ \ \textit{
 Let $M$ be the cellular model of an orientable relative presentation $(C_n,x:W)$ for a group $E$, where $C_n$ is a cyclic group of order $n$. Assume that $\pi_2M = 0$. If $\nu^f:E \ra C_n$ is any retraction onto the coefficient group $C_n$, then $C_n$ acts freely via the shift on the nonidentity elements of the cyclically presented group $G_n(\rho^f(W))$.}

\vspace{.25in}

\noindent The proof of Theorem A depends on a relation between the homotopy condition $\pi_2M = 0$ and combinatorial asphericity of the cyclic presentation $\mathcal{P}_n(\rho^f(W))$. The following result extends \cite[Lemma 3.1]{GH} and \cite[Corollary 4.3]{EdjvetWilliams}.

\vspace{.25in}

\noindent \textbf{Theorem \ref{asphTransfer}}\ \ \textit{Let $M$ be the cellular model of an orientable relative presentation $(C_n,x:W)$ for a group $E$ where $C_n$ is the cyclic group of order $n$ generated by $a$ and
$$
W = x^{\epsilon_1}a^{p_1}\ldots x^{\epsilon_L}a^{p_L} \in C_n \ast \<x\>
$$
is cyclically reduced. Suppose that $\nu^f:E \ra C_n$ is a retraction given by $\nu^f(a) = a$ and $\nu^f(x) = a^f$. Then, $\rho^f(W)$ is cyclically reduced and the cyclic presentation $\mathcal{P}_n(\rho^f(W))$ is orientable. Moreover, $\mathcal{P}_n(\rho^f(W))$ is combinatorially aspherical if and only if $\pi_2M = 0$.}

\vspace{.25in}

Even for orientable cyclic presentations, the converse to Theorem A is not true. For example, the ``original" Fibonacci group \cite{JC} was $G_5(x_0x_1x_2^{-1})$, which turned out to be cyclic of order $11$. The associated cyclic presentation is not combinatorially aspherical, but a calculation with GAP \cite{GAP} shows that the orbits of the shift partition the nonidentity elements of $G_5(x_0x_1x_2^{-1})$ into two disjoint $5$-cycles.

The converse to Theorem A does hold for a class of cyclically presented groups that was studied by M.~Edjvet and G.~Williams \cite{EdjvetWilliams}. Given integers $k$ and $l$ modulo a positive integer $n$, let $G_n(k,l)$ be the group given by the (orientable) cyclic presentation

$$
\mathcal{P}_n(k,l) = \mathcal{P}_n(x_0x_kx_l)
$$
that has generators $x_0,\ldots,x_{n-1}$ and relators of the form $x_ix_{i+k}x_{i+l}$ for $i = 0,\ldots,n-1$.

\vspace{.25in}

\noindent \textbf{Theorem B}\ \ \textit{The cyclic group $C_n$ acts freely via the shift on the nonidentity elements of $G_n(k,l)$ if and only if the presentation $\mathcal{P}_n(k,l)$ is combinatorially aspherical.}

\vspace{.25in}

\noindent \textbf{Theorem C}\ \  \textit{The group $G_n(k,l)$ is finite if and only if the shift $\theta$ has a nonidentity fixed point.}

\vspace{.25in}

\noindent Edjvet and Williams \cite{EdjvetWilliams} developed a taxonomy for the parameters $(n,k,l)$ that characterizes finiteness for the groups $G_n(k,l)$ and (topological) asphericity for the two-dimensional cellular models of the presentations $\mathcal{P}_n(k,l)$. (See also \cite{CRS}.) The proofs of Theorems B and C follow this taxonomy to establish the stated connections to the dynamics of the shift.

There are variations in the use of the term ``aspherical" and we navigate these rather carefully. The results presented here are consistent with those of \cite{EdjvetWilliams}, where the parameters $k$ and $l$ are restricted to be nonzero. I am grateful to Gerald Williams for his suggestions about this work.

\section{Asphericity}\label{secAsph}\label{sec:Asphericity}

The following topological preliminaries can be found in \cite{AJSAlgTop}. A topological space $Y$ is \textbf{aspherical} if each spherical map $S^k \ra Y$, $k \geq 2$, is homotopic to a constant map. A connected aspherical CW complex with fundamental group $G$ is a \textbf{$K(G,1)$-complex}; its homotopy type is uniquely determined by $G$. The \textbf{cellular model} of a group presentation $(\mathbf{x}:\mathbf{r})$ is a connected two-dimensional CW complex $K$ whose one-skeleton $K^{(1)}$ consists of a single  zero-cell $\ast$ together with a single one-cell $c^1_x$ for each generator $x \in \bf{x}$. Selection of an orientation for each one-cell determines an isomorphism $\pi_1K^{(1)} \cong F(\mathbf{x})$ of the fundamental group of the one-skeleton $K^{(1)}$ with the free group $F = F({\bf x})$ with basis $\mathbf{x}$. The cellular model $K$ has a two-cell $c^2_r$ for each relator $r \in \mathbf{r}$ whose boundary attaching map spells the relator $r \in F$. The fundamental group of the cellular model is isomorphic to the group given by the presentation.

%The following topological preliminaries can be found in \cite{AJSAlgTop}. A topological space $Y$ is \textbf{aspherical} if each spherical map $S^k \ra Y$, $k \geq 2$, is homotopic to a constant map. The homotopy type of a connected aspherical CW complex is uniquely determined by the isomorphism type of its fundamental group. The \textbf{cellular model} of a group presentation $(\mathbf{x}:\mathbf{r})$ is a connected two-dimensional CW complex $K$ whose one-skeleton $K^{(1)}$ consists of a single  zero-cell $\ast$ together with a single one-cell $c^1_x$ for each generator $x \in \bf{x}$. Selection of an orientation for each one-cell determines an isomorphism $\pi_1(K^{(1)},\ast) \cong F(\mathbf{x})$ of the fundamental group of the one-skeleton $K^{(1)}$ with the free group with basis $\mathbf{x}$. The cellular model $K$ has a two-cell $c^2_r$ for each relator $r \in \mathbf{r}$ with characteristic map $\phi_r:B^2 \ra K$ whose boundary attaching map $\dot{\phi}_r:S^1 \ra K^{(1)}$ spells the relator $r \in F(\mathbf{x}) \equiv \pi_1(K^{(1)},\ast)$. The inclusion-induced homomorphism $\pi_1(K^{(1)},\ast) \ra \pi_1(K,\ast)$ is surjective with kernel equal to the smallest normal subgroup of $F(\mathbf{x})$ that contains $\mathbf{r}$. Thus the fundamental group of the cellular model is isomorphic to the group given by the presentation.

\paragraph{Combinatorial Asphericity:} Identity sequences \cite{CCH,AJSFrame} describe spherical maps into cellular models of group presentations. Given a presentation $(\mathbf{x}:\mathbf{r})$ for a group $G$, the free group $F = F(\mathbf{x})$ with basis $\mathbf{x}$ acts on the free group $\mathbb{F} = \mathbb{F}(\mathbf{x}:\mathbf{r}) = F(F\times \mathbf{r})$ with basis $F \times \mathbf{r}$ by $v \cdot (u,r) = (vu,r)$. The homomorphism $\partial: \mathbb{F} \ra F$ given by $\partial(u,r) = uru^{-1}$ is $F$-equivariant where $F$ acts on itself by (left) conjugation. The image of $\partial$ is the normal closure of $\mathbf{r}$ in $F$ so the cokernel of $\partial$ is isomorphic to $G$. Elements of the kernel $\mathbb{I} = \mathbb{I}(\mathbf{x}:\mathbf{r}) = \ker \partial$ are called \textbf{identity sequences} for the presentation $(\mathbf{x}:\mathbf{r})$ \cite{CCH}. Given $u,v \in F$, $r,s \in \mathbf{r}$, and $\epsilon,\delta = \pm 1$, there is the Peiffer identity $(u,r)^\epsilon(v,s)^\delta(u,r)^{-\epsilon}(ur^\epsilon u^{-1}v,s)^{-\delta} \in \mathbb{I}(\mathbf{x}:\mathbf{r})$. Let $\mathbb{P} = \mathbb{P}(\mathbf{x}:\mathbf{r})$ denote the normal closure in $\mathbb{F}$ of the set of all Peiffer identities. The $F$-action on $\mathbb{F}$ descends to a $G$-action on the abelian quotient group $\mathbb{I}/\mathbb{P}$. Here is a working definition of combinatorial asphericity:

\begin{Proposition}\emph{\textbf{(\cite[Proposition 1.4]{CCH})}}\label{CA}
  A presentation $(\mathbf{x}:\mathbf{r})$ for a group $G$ is \textbf{combinatorially aspherical} if and only if $\mathbb{I}(\mathbf{x}:\mathbf{r})/\mathbb{P}(\mathbf{x}:\mathbf{r})$ is generated as a $\mathbb{Z}G$-module by (classes of) length two identity sequences, that is by identity sequences of the form

  $$
  (u,r)^\epsilon(v,s)^\delta
  $$
  where $u,v \in F$, $r,s \in \bf{r}$, and $\epsilon,\delta = \pm 1$.
\end{Proposition}

%Sieradski \cite[Section 2]{AJSFrame} (see also \cite{AJSAlgTop}) describes an isomorphism $\pi_2 K \cong \mathbb{I}/\mathbb{P}$ of $\mathbb{Z}G$-modules where $K$ is the two-complex modeled on the presentation $(\mathbf{x}:\mathbf{r})$.
%
%Letting $K = \bigvee_\mathbf{x}S^1_x \cup \bigcup_\mathbf{r} c^2_r$ be the cellular model of the presentation $(\mathbf{x}:\mathbf{r})$, there is an $F$-equivariant homomorphism $\eta: \mathbb{F}(\mathbf{x}:\mathbf{r}) \ra \pi_2(K,K^{(1)},\ast)$ that maps the basis element $(1,r) \in \mathbb{F}(\mathbf{x}:\mathbf{r})$, $r \in \mathbf{r}$, to the class $[\phi_r:(B^2,S^1)\ra(K,K^{(1)})]$ of the characteristic map for the two-cell $c^2_r$ in the relative homotopy group $\pi_2(K,K^{(1)},\ast)$ with its homotopy action by $\pi_1(K^{(1)},\ast) \equiv F$. J.~H.~C.~Whitehead showed \cite{JHC} that $\eta$ is surjective with kernel equal to the group $\mathbb{P}$ of Peiffer identities. See also \cite[Lemma 2.5]{AJSAlgTop} or \cite[page 127]{AJSFrame}. Using the long exact homotopy sequence for the based pair $(K,K^{(1)},\ast)$, this leads to an isomorphism $\pi_2 (K,\ast) \cong \mathbb{I}/\mathbb{P}$ of $\mathbb{Z}G$-modules. For an identity sequence $\sigma \in \mathbb{I}$, we denote the corresponding homotopy class by $[\sigma] \in \pi_2 (K,\ast)$.

Given a presentation $(\mathbf{x}:\mathbf{r})$ for a group $G$ and having cellular model $K$, there is an isomorphism

$$
\pi_2 K \cong \mathbb{I}(\mathbf{x}:\mathbf{r})/\mathbb{P}(\mathbf{x}:\mathbf{r})
$$
of $\mathbb{Z}G$-modules. (See \cite[Section 2]{AJSFrame} or \cite{AJSAlgTop}, where this is attributed to K.~Reidemeister \cite{R}.) We sometimes write $[\sigma] = \sigma\mathbb{P} \in \mathbb{I}/\mathbb{P}$. Since $K$ is two-dimensional, it follows that $K$ is aspherical (in the topological sense) if and only if $\mathbb{I}(\mathbf{x}:\mathbf{r}) = \mathbb{P}(\mathbf{x}:\mathbf{r})$. This implies that a presentation with aspherical cellular model is combinatorially aspherical. However, the cellular model of a combinatorially aspherical presentation can support essential spheres.

A relator $r \in \mathbf{r}$ is a \textbf{proper power} if it is freely equal to $r = \dot{r}^{e}$ where $e > 1$. When the exponent $e$ is maximal it is the \textbf{exponent} of $r$ and $\dot{r}$ is the \textbf{root}. A proper power relator gives rise to an identity sequence $(1,r)(\dot{r},r)^{-1}$.

\begin{Example}\label{Dn3} It is well known and fairly easy to show that the presentation $(a:a^n)$ for the cyclic group $C_n$ of order $n$ is combinatorially aspherical by showing that

$$
\mathbb{I}(a:a^n)/\mathbb{P}(a:a^n) = \mathbb{Z}C_n \cdot [\Delta_n]\ \ \mbox{where}\ \ \Delta_n = (1,a)(a,a^n)^{-1}.
$$
A consequence is that a $K(C_n,1)$-complex $D_n$ can be constructed that has just a single three-cell with attaching map $S^2 \ra D_n^{(2)}$ determined by the identity sequence $\Delta_n \in \mathbb{I}(a:a^n)$. \qed
\end{Example}

A relator of $(\mathbf{x}:\mathbf{r})$ is \textbf{redundant} if it is freely trivial or else is freely conjugate to another relator or its inverse. No relator in a combinatorially aspherical presentation can be freely trivial (as every identity sequence has even length), but the concept does allow for the possibility that one relator $r$ is freely conjugate to another relator $s$ or its inverse: $r = vs^\delta v^{-1}$. This leads to identity sequences of the form $(1,r)(v,s)^{-\delta}$. When $r$ and $s$ are cyclically reduced, this simply means that $r$ is a cyclic permutation of $s^\delta$ (and $v$ is an initial segment of $r$). Identity sequences arising from proper power or redundant relators are called Peiffer identities of the second kind in \cite{CCH,Hueb}.

To illustrate the possibilities, the cyclic presentation

$$
\mathcal{P}_3(1,2) = \mathcal{P}_3(x_0x_1x_2) = (x_0,x_1,x_2:x_0x_1x_2,x_1x_2x_0,x_2x_0x_1)
$$
for the free group of rank two has two redundant relators. In the cyclic presentation $\mathcal{P}_n(0,0) = \mathcal{P}_n(x_0^3)$ for the free product of $n$ copies of the cyclic group of order three, each relator has exponent three and root of length one. It is not difficult to show that the presentations $\mathcal{P}_3(x_0x_1x_2)$ and $\mathcal{P}_n(x_0^3)$ are combinatorially aspherical. %It is also possible to derive both claims simultaneously from considerations involving aspherical relative presentations.

\paragraph{Aspherical Relative Presentations:} Consider a group presentation of the form

$$
\mathcal{Q} = (a,x:a^n,W(a,x)).
$$
If we view the word $W(a,x)$ as an element $W$ in the free product $C_n \ast \<x\>$ of $C_n$ with the infinite cyclic group generated by $x$, we obtain a \textbf{relative presentation}

%\begin{equation}\label{relPres}
%  \mathcal{R} = (C_n,x:W)
%\end{equation}
$$
\mathcal{R} = (C_n,x:W)
$$
with \textbf{coefficient group} $C_n$, in the sense of \cite{BP}. Conversely, the ordinary presentation $\mathcal{Q}$ is a \textbf{lift} of $\mathcal{R}$. Both $\mathcal{Q}$ and $\mathcal{R}$ define the same group:

$$
E = F(a,x)/\<\<a^n, W(a,x)\>\> \cong (C_n\ast\<x\>)/\<\<W\>\>.
$$
The one-relator relative presentation $\mathcal{R}$ is \textbf{orientable} if $W$ is not conjugate to $W^{-1}$ in $C_n \ast \<x\>$. Orientability means that $W$ is not conjugate to a word of the form $X\alpha X^{-1}\beta$ where $X \in C_n \ast \<x\>$ and $\alpha^2=\beta^2=1$ in $C_n$. See \cite[(1.3)]{BP}.

A cellular model for the relative presentation $\mathcal{R}$ was described in \cite[Section 4]{BP}. Working in the free product $C_n \ast \<x\>$, first write $W = \dot{W}^e$ where $\dot{W} \in C_n \ast \<x\>$ and $e$ is maximal, so again $\dot{W}$ is the \textbf{root} and $e$ is the \textbf{exponent} of $W$. Let $D_n$ be a $K(C_n,1)$-complex with two-skeleton modeled on the presentation $(a:a^n)$ and having just a single three-cell as in Example \ref{Dn3}. In the same way, let $D_e$ be a $K(C_e,1)$-complex, unless $e=1$, in which case $D_e$ is a disc. Form the one point union $D_n \vee S^1_x$ with fundamental group $\pi_1(D_n \vee S^1_x) \cong C_n \ast \<x\>$. Now the one-skeleton $D_e^{(1)} = S^1_e$ is a circle and there is a loop $\omega: S^1_e \ra D_n \vee S^1_x$ that reads the root $\dot{W} \in C_n \ast \<x\>$. The \textbf{cellular model} of $\mathcal{R}$ is the pushout $M$ obtained from the union of $D_n \vee S^1_x$ and $D_e$ by identifying points of $D_e^{(1)} = S^1_e$ with their images in $D_n \vee S^1_x$ under $\omega$:

$$
M = \left(D_n \vee S^1_x\right) \cup_{z \sim \omega(z)} D_e.
$$
%Note that $M$ is obtained from the complex $X$ in the proof of Theorem \ref{splitCrit} by adjoining the cells of of $D_e$ in dimensions three and higher: $X = D_n \cup M^{(2)}$. In particular, if $W$ is not a proper power in $C_n \ast \<x\>$, so that $e=1$, then  $X = M$.

%The relative presentation $\mathcal{R}$ is \textbf{aspherical} if the relative homotopy group $\pi_2(M,D_n)$ is trivial. (See \cite[Theorem 4.1]{BP}.) Since $D_n$ is aspherical, the exact homotopy sequence for the pair $(M,D_n)$ reveals that $\mathcal{R}$ is aspherical if and only if the element $a$ has order exactly $n$ in $E$ and $\pi_2 M = 0$. In fact, the cellular model $M$ is (toplogically) aspherical and this leads to the main results on for groups given by aspherical orientable relative presentations.

\begin{Proposition}\label{relComb} Let $M$ be the cellular model of a relative presentation $\mathcal{R} = (C_n,x:W)$ for a group $E$. Let $\mathcal{Q} = (a,x:a^n,W(a,x))$ be a lift of $\mathcal{R}$ where $W(a,x)$ is obtained from $W \in C_n \ast \<x\>$ by replacing each appearance in $W$ of a coefficient from $C_n$ by the appropriate  subword of the form $a^k$ where $0 \leq k < n$.
\begin{enumerate}
  \item[(a)] If $\pi_2M = 0$, then $\mathcal{Q}$ is combinatorially aspherical.
  \item[(b)] If $\mathcal{R}$ is orientable and $\mathcal{Q}$ is combinatorially aspherical, then $\pi_2M = 0$.
  \item[(c)] If $\mathcal{R}$ is orientable and $\pi_2M = 0$, then $M$ is aspherical: $\pi_k M = 0$, $k\geq 2$.
\end{enumerate}
\end{Proposition}

\proof The careful way that $W(a,x)$ is lifted from $W$ ensures that the exponent of $W$ in $C_n \ast \<x\>$ is the same as the exponent of $W(a,x)$ in the free group on $a$ and $x$. The two-skeleton $M^{(2)}$ is the cellular model of the ordinary presentation $\mathcal{Q}$ and the three-cells of $M$ are attached by spherical maps $S^2 \ra M^{(2)}$ that correspond to the identity sequences

$$
(1,a^n)(a,a^n)\ \ \mbox{and}\ \ (1,W(a,x))(\dot{W}(a,x),W(a,x))^{-1}
$$
for the presentation $\mathcal{Q}$. If $\pi_2M = 0$, as in (a), then the cellular approximation theorem implies that $\pi_2 M^{(2)}$ is generated as a $\mathbf{Z}E$-module by the homotopy classes of these attaching maps and so $\mathcal{Q}$ is combinatorially aspherical by Proposition \ref{CA}. The hypotheses in (b) imply that $\pi_2M^{(2)}$ is $\mathbb{Z}E$-generated by classes of spherical maps corresponding to identity sequences over $\mathcal{Q}$ having the form $(1,a^n)(a^k,a^n)^{-1}$ and $(1,W(a,x))(\dot{W}(a,x)^k,W(a,x))^{-1}$. This uses the fact that the centralizers of $a^n$ and $W(a,x)$ in the free group with basis $a$ and $x$ are generated by $a$ and $\dot{W}(a,x)$ respectively. However these spherical maps are inessential in $D_n$ and in (the image in $M$ of) $D_e$, respectively, and hence in $M$. Thus $\pi_2M = 0$. The statement (c) is a direct consequence of \cite[Theorem 4.2]{BP}, which used \cite[Theorem 4.2]{HowieHighPower}. \qed

%For (b), assume that $\mathcal{R}$ is orientable and $\pi_2M = 0$. By (a), the homomorphisms induced by the inclusions of $D_n \vee S^1_x$ and $D_e$ into $M = (D_n \vee S^1_x) \cup D_e$ are both free groups. Since $\pi_2M = 0$, a theorem of J.~Howie \cite[Theorem 4.2]{HowieHighPower} implies that $M$ is aspherical.
%
%\vspace{.25in}
%
%The following results are from \cite{BP}.

\begin{Theorem}[\cite{BP}] \label{OrAsphMain} Let $\mathcal{R} = (C_n,x:W)$ be an orientable relative presentation for a group $E$. Assume that $\pi_2M = 0$ where $M$ is the cellular model $M$ of $\mathcal{R}$.
\begin{enumerate}
  \item[(a)] The root $\dot{W}$ generates a subgroup $C_e$ of $E$ with order equal to the exponent of $W$ in $C_n \ast \<x\>$; the coefficient group $C_n$ embeds naturally in $E$.
  \item[(b)] Given any left $ZE$-module $N$ and any $k \geq 3$, the canonical map $H^k(E;N) \ra H^k(C_n;N) \times H^k(C_e;N)$ is an isomorphism.
  \item[(c)] Any finite subgroup of $E$ is conjugate to a subgroup of either $C_n$ or $C_e$.
  \item[(d)] All conjugates of $C_n$ have trivial intersection with all conjugates of $C_e$. If $g \in E$ and $gC_ng^{-1} \cap C_n \neq 1$, then $g \in C_n$. If $g \in E$ and $gC_eg^{-1} \cap C_e \neq 1$, then $g \in C_e$.
\end{enumerate}
\end{Theorem}

\proof These results were all proved in \cite{BP}. The result (a) follows from Proposition \ref{relComb}(a) by \cite[Proposition 1]{Hueb}. The calculation (b) follows from Proposition \ref{relComb}(c) by standard methods because the cells of the cellular model $M$ in dimensions three and up are precisely those of $D_n$ and $D_e$. Characterization of the finite subgroups as in (c) follows from the cohomology calculation (b) and a theorem of J.-P.~Serre that appeared in \cite{Hueb}. The conjugacy results in (d) can be extracted either from Serre's theorem or from \cite[Theorem 5]{HS}. A nice alternative approach to (c) and (d) for the case where $W$ is not a proper power is given in \cite{FR}. \qed

\vspace{.25in}

The concept of asphericity for relative presentations was formulated in \cite{BP} in terms of a geometric condition involving pictures over relative presentations. For orientable relative presentations, that geometric condition implies the homotopy condition $\pi_2M = 0$ \cite[Theorem 4.2]{BP}. This is not the case for nonorientable relative presentations.

\begin{Example}\label{ex:nonor}
  If $1 \neq a \in C_2$, then the relative presentation $\mathcal{R} = (C_2,x:[x,a])$ for the direct product $E = C_2 \times \<x\>$ is aspherical in the sense of \cite{BP}, but is not orientable. Here $[x,a]$ denotes the commutator $xax^{-1}a^{-1}$. The group $C_2 \times \<x\>$ admits a conjugacy relation $xC_2x^{-1} = C_2$, so the conclusions of Theorem \ref{OrAsphMain}(b,d) do not extend to nonorientable relative presentations that are aspherical in the sense of \cite{BP}. Moreover, it follows from \cite[Theorem 3]{Hueb} that the ordinary presentation $\mathcal{Q} = \left(a,x:a^2,xax^{-1}a^{-1}\right)$ is not combinatorially aspherical. Indeed, $\mathcal{Q}$ supports an identity sequence

$$
(1,a^2)(x,a^2)^{-1}(1,[x,a])(a,[x,a])
$$
that arises from the free commutator identity $[x,a^2] = [x,a]a[x,a]a^{-1}$ and which represents a nontrivial homotopy class in $\pi_2M$.
\end{Example}

\section{Retractions}\label{sec:Retrations}

Suppose that $\nu: E\ra C$ is a retraction of a group $E$ onto a subgroup $C$. Thus $\nu(c) = c$ for all $c \in C$. Let $G$ denote the kernel of $\nu$ so that $E \cong G \rtimes C$. Now $C$ acts on $G$ by (left) conjugation (that is, $c \cdot g = cgc^{-1}$) and also on the coset space $E/C$ by left multiplication (that is, $c \cdot wC = cwC$).

\begin{Lemma}\label{TransferDynamic}
  $E/C \cong G$ as left $C$-sets with $1C \in E/C$ corresponding to $1 \in G$.
\end{Lemma}

\proof The desired isomorphism is given by the $C$-equivariant function $E/C \ra G$ that maps $wC \in E/C$ to $w \cdot \nu(w)^{-1} \in G$ for each $w \in E$. \qed

\vspace{.25in}

Given a cyclic presentation $\mathcal{P}_n(w)$, the semidirect product $E_n(w) = G_n(w) \rtimes_\theta C_n$ has a two-generator two-relator presentation of the form $\mathcal{Q} = (x,a:a^n, W(a,x))$ where $W(a,x)$ is obtained from $w(x_0,\dots,x_{n-1})$ by rewriting $x_0 = x$ and $x_i = a^ixa^{-i}$ for $i=1,\dots,n-1$. This is because the rewrites of the shifts $\theta^v(w)$ are pairwise conjugate in the free group on $a$ and $x$ and so are redundant. There is a retraction $\nu = \nu^0: E_n(w) \ra C_n$ given by $\nu(a) = a$ and $\nu(x) = a^0 = 1$. There can be other retractions, as follows.

\begin{Lemma}\label{Lemma:splitCrit}
  Let $\mathcal{Q} = \left(a,x:a^n, W(a,x)\right)$ be a presentation for a group $E$. Let $p$ and $\epsilon$ be the exponent sums of $a$ and $x$ in $W(a,x)$, respectively. The assignments

  $$\
  \nu^f(a) = a\ \ \mbox{and}\ \ \nu^f(x) = a^f
  $$
  define a retraction of $E$ onto the cyclic subgroup $C_n$ of order $n$ if and only if $\epsilon f + p$ is divisible by $n$. \qed
\end{Lemma}

\vspace{.25in}

Cyclic presentations for retraction kernels arise from the following Reidemeister-Schreier rewriting process \cite[Theorem 2.9, page 94]{MKS} for words of the form

\begin{equation}\label{Wform}
  W(a,x) = x^{\epsilon_1}a^{p_1}\ldots x^{\epsilon_L}a^{p_L}
\end{equation}
where $L \geq 1$, $\epsilon_i = \pm 1$, and the $p_i$ are integers. Given an integer $f$, the process involves two integer-valued functions. The first function is defined by $v(1) = 0$ and

$$
v(i) = \sum_{j=0}^{i-1} \epsilon_{j}f + p_{j}.
$$
for $i = 2,\dots,L$.
%$$
%v(i+1) = v(i) + \epsilon_{i}f + p_{i}.
%$$
The second function is defined in terms of the first by

$$
u(i) = v(i) + \frac{\epsilon_i - 1}{2}f
= \left\{ \begin{array}{ll}
            v(i) & \mbox{if}\ \ \epsilon_i = 1\\
            v(i)-f & \mbox{if}\ \ \epsilon_i = -1.
            \end{array}
            \right.
$$
Given the integer $f$ and a word $W(a,x)$, the rewritten word $\rho^f(W(a,x))$ is defined to be

$$
\rho^f(W(a,x)) = x_{u(1)}^{\epsilon_1}\ldots x_{u(L)}^{\epsilon_L}.
$$
%sLThe rewritten word is defined on the alphabet $X = \{x_0,\ldots,x_{n-1}\}$.
%Interpreting subscripts modulo $n$, the rewriting process respects the shift $\theta(x_i) = x_{i+1}$ in the sense that
%
%$$
%\rho^f(W) = \theta(\rho^f(W)).
%$$
Some examples:

\begin{eqnarray*}
% \nonumber to remove numbering (before each equation)
  \rho^f(xa^pxa^qxa^r) &=& x_0x_{f+p}x_{2f+p+q}, \\
   \rho^f(x^2a^pxa^q) &=& x_0x_fx_{2f+p}, \\
  \rho^f(x^2a^px^{-1}a^q) &=& x_0x_fx^{-1}_{f+p},\\
  \rho^{0}(x^3) &=& x_0^3,\\
  \rho^1(x^3) &=& x_0x_1x_2.
\end{eqnarray*}
\noindent When a positive integer $n$ is fixed, $W(a,x)$ determines $W \in C_n \ast \<x\>$. By interpreting subscripts modulo $n$, the rewriting process yields a word $\rho^f(W)$ on the alphabet $\{x_0,\ldots,x_{n-1}\}$. In its turn, the cyclic presentation $\mathcal{P}_n(\rho^f(W))$ is defined and depends on $n$. Thus $\mathcal{P}_6(\rho^2(x^3)) = \mathcal{P}_6(x_0x_2x_4) = \mathcal{P}_6(2,4)$, whereas $\mathcal{P}_3(\rho^2(x^3)) = \mathcal{P}_3(x_0x_2x_1) = \mathcal{P}_3(2,1)$.

%\begin{eqnarray*}
%% \nonumber to remove numbering (before each equation)
%  \mathcal{P}_6(\rho^2((xa)^3)) &=& \mathcal{P}_6(x_0x_3x_0) = \mathcal{P}_6(3,0), \\
%  \mathcal{P}_3(\rho^1((xa)^3)) &=& \mathcal{P}_3(x_0x_2x_1) = \mathcal{P}_3(2,1).
%\end{eqnarray*}
%

%The Lemma applies with $n = \infty$, in which case the relator $a^n$ is omitted and the result is about the one-relator group $E = (a,x : W)$. The only integer that is divisible by $n = \infty$ is zero. The cyclic group $C_n = C_\infty$ is infinite cyclic. The cyclic presentation $\mathcal{P}_\infty(\rho_1(W(a,x)))$ uses the infinite generating set $\widehat{X}$ and the split extension $E \cong G_\infty(\rho_1(W(a,x))) \rtimes_{\widehat{\theta}} C_\infty$ involves the infinite order shift $\widehat{\theta}$.\\

\begin{Theorem}\label{splitCrit} Let $(C_n,x:W)$ be a relative presentation for a group $E$ where $C_n$ is a cyclic group of order $n$ generated by $a$. If $\nu^f:E \ra C_n$ is a retraction given by $\nu^f(a) = a$ and $\nu^f(x) = a^f$, then the element $x_0 = xa^{-f}$ normally generates $\ker \nu^f$ in $E$ and the conjugates $x_i = a^{i}x_0a^{-i}$ ($0 \leq i \leq n-1$) form a generating set in a cyclic presentation $\mathcal{P}_n(\rho^f(W))$ for $\ker \nu^f$. In particular,

$$
E \cong G_n(\rho^f(W)) \rtimes_\theta C_n
$$
and $\theta(g) = aga^{-1}$ in $E$ for all $g \in \ker \nu^f$.
\end{Theorem}

\proof Let $M$ be the cellular model of the relative presentation $(C_n,x:W)$ with a single zero-cell $\ast \in M^{(0)}$ so that $\pi_1(M,\ast) \cong E$. There is regular covering projection $p:\widehat{M} \ra M$ and a zero-cell $\ast_0 \in p^{-1}(\ast) \subseteq \widehat{M}$ such that $p_\sharp(\pi_1(\widehat{M},\ast_0)) = \ker \nu^f \unlhd \pi_1(X,\ast) = E$ where $p_\sharp:\pi_1(\widehat{M},\ast_0) \ra \pi_1(M,\ast)$ is the homomorphism of fundamental groups induced by $p$. By lifting paths through the covering projection $p$, the choice of basepoint $\ast_0$ and the homomorphism $\nu^f$ determine an action of $E$ on $\widehat{M}$ that identifies $\nu^f(E) = C_n$ as the automorphism group $\mathrm{Aut}(p)$ of the covering projection $p$. %The $E$-action on $\widehat{M}$ is given by path lifting.

For each cell of the base complex $M$, the automorphism group $\mathrm{Aut}(p) \cong C_n$ freely permutes a single orbit of cells lying over it in $\widehat{M}$. Thus the covering complex $\widehat{M}$ has zero-cells $\ast_i = a^i \cdot \ast_0$ for $i = 0,\ldots, n-1$, which constitute the full pre-image $p^{-1}(\ast)$ lying over the zero-cell $\ast$ of $M$. In the one-skeleton $\widehat{M}^{(1)}$, there are one-cells $a_0,\ldots,a_{n-1}$ lying over the one-cell $a$ of $M$ and one-cells $x_0,\ldots,x_{n-1}$ lying over $x$. The endpoints of the one-cells are determined by the homomorphism $\nu^f$ and in turn reflect the action of $E = \<a,x\>$ on $\widehat{M}$. Thus the one-cell $a_i$ has initial zero-cell $\ast_i$ and terminal zero-cell $a\cdot \ast_i = \ast_{i+1}$. The one-cell $x_i$ has initial zero-cell $\ast_i$ and terminal zero-cell $x \cdot \ast_i = a^f \cdot \ast_i = \ast_{i+f}$.

Since $\nu^f$ maps the subgroup $\pi_1(D_n,\ast) = \<a\> \cong C_{n}$ isomorphically onto $C_{n} \cong \mathrm{Aut}(p)$, the pre-image $p^{-1}(D_n)$ is isomorphic to the universal covering complex $\widetilde{D}_n$ of $D_n$. The quotient map that collapses the contractible $C_n$-equivariant subcomplex $\widetilde{D}_n$ of $\widehat{M}$ to a point is a homotopy equivalence onto the quotient complex $K = \widehat{M}/\widetilde{D}_n$, which has a single zero-cell. The action of $C_n$ on $\widehat{M}$ descends to a cellular $C_n$-action on $K$ that is free away from the zero-cell. This means that $K$ models a cyclic presentation for $\pi_1(\widehat{M},\ast_0) \cong \ker \nu^f$. In the quotient complex $K$, there is a single zero-cell, together with a free $C_n$-orbit of one-cells, still denoted $x_0, \ldots, x_{n-1}$. The two-cells of $p^{-1}(c^2_W) \subseteq \widehat{M}$ descend to a free $C_n$-orbit of two-cells in $K$ with boundary paths $\theta^v(\rho^f(W))$ where $v = 0,\ldots,n-1$. Thus the two-skeleton $K^{(2)}$ is the two-complex modeled on the presentation $\mathcal{P}_n(\rho^f(W))$, as claimed. Viewed inside the group $E$ with generators $a$ and $x$, the generator $x_i$ for $\ker \nu^f = G_n(\rho^f(W))$ is given by $x_i = a^{i}xa^{-(i+f)}$, which describes a loop in $M^{(1)}$ that lifts to a loop in $\widehat{M}^{(1)}$ based at $\ast_0$ that in turn lies over the one-cell $x_i$ of the quotient complex $K$. Thus $\ker\nu^f$ is normally generated in $E$ by $x_0 = xa^{-f}$ and the shift $\theta$ for the cyclically presented group $\ker \nu^f = G_n(\rho^f(W))$ is given by $\theta(g) = aga^{-1}$ for all $g \in \ker \nu^f$. \qed

\vspace{.25in}

There are many results on asphericity of relative presentations \cite{BBP,BP,D,ECap,HM,Met} that can potentially be applied to the study of the shift on cyclically presented groups.

\begin{Theorem}\label{freeAction} Let $M$ be the cellular model of an orientable relative presentation $(C_n,x:W)$ for a group $E$, where $C_n$ is a cyclic group of order $n$. Assume that $\pi_2M = 0$. If $\nu^f:E \ra C_n$ is any retraction onto the coefficient group $C_n$, then $C_n$ acts freely via the shift on the nonidentity elements of the cyclically presented group $G_n(\rho^f(W))$.
\end{Theorem}

\proof Suppose $g \in G_n(\rho^f(W))$ and $\theta^k(g) = g$ where $1 \leq k < n$. By Theorem \ref{splitCrit}, the shift $\theta$ on $G_n(\rho^f(W))$ is realized as conjugation by $a \in C_n$ in $E \cong G_n(\rho^f(W)) \rtimes_\theta C_n$. Thus $a^kga^{-k} = g$ and hence $1 \neq a^k = ga^kg^{-1} \in C_n \cap gC_ng^{-1}$. By Theorem \ref{OrAsphMain}(c), $g \in C_n \cap G_n(\rho^f(W)) = 1$. \qed

\section{Proof of Theorem A}

A cyclic presentation $\mathcal{P}_n(w)$ is orientable unless $w$ is freely equal to the inverse of one its shifts.

\begin{Example}
   As a companion to the nonorientable relative presentation in Example \ref{ex:nonor}, consider $\rho^0(xax^{-1}a^{-1}) = x_0x_1^{-1}$ and the nonorientable cyclic presentation

$$
\mathcal{P}_2(x_0x_1^{-1}) = (x_0,x_1:x_0x_1^{-1}, x_1x_0^{-1})
$$
for the infinite cyclic group, which has two mutually inverse defining relators. This presentation is combinatorially aspherical but the shift acts as the identity. For nonorientable cyclic presentations, the shift almost never acts freely on the nonidentity elements. \qed
\end{Example}

\begin{Lemma}\label{orientableShape}
  Let $w = w(x_0,\ldots, x_{n-1})$ be a nonempty cyclically reduced word. The cyclic presentation $\mathcal{P}_n(w)$ is nonorientable if and only if $n=2m$ is even and $w = w(x_0,\ldots,x_{n-1})$ is freely equal to $u\theta^m(u)^{-1}$ for some linearly reduced word $u = u(x_0,\ldots,x_{n-1})$.
\end{Lemma}

\proof If $n=2m$ and $w = u\theta^m(u)^{-1}$ then $w\theta^m(w)$ is freely trivial so $\mathcal{P}_n(w)$ is not orientable. Conversely, suppose that $w\theta^{v}(w)$ is freely trivial where $1\leq v < n$. Write $w = x_{p_1}^{\epsilon_1}\ldots x_{p_L}^{\epsilon_L}$, so that

$$
x_{p_1}^{\epsilon_1}\ldots x_{p_L}^{\epsilon_L}x_{p_1+v}^{\epsilon_1}\ldots x_{p_L+v}^{\epsilon_L}=1.
$$
Since $w$ is reduced, cancellation implies that $\epsilon_{L-j} = -\epsilon_{j+1}$ and $p_{L-j} \equiv p_{j+1}+v$ modulo $n$ for $j = 0,\ldots,L-1$. The conditions on the $\epsilon$'s ensure that $L=2M$ is even. The conditions on the $p$'s show that $p_L \equiv p_1+v \equiv p_L+2v$ modulo $n$, which means that $n = 2v$ is even. Letting $u = x_{p_1}^{\epsilon_1}\ldots x_{p_M}^{\epsilon_M}$, which is linearly reduced, we have

$$
u\theta^v(u)^{-1} = x_{p_1}^{\epsilon_1}\ldots x_{p_M}^{\epsilon_M}\left(x_{p_1+v}^{\epsilon_1}\ldots x_{p_M+v}^{\epsilon_M}\right)^{-1} = x_{p_1}^{\epsilon_1}\ldots x_{p_M}^{\epsilon_M}x_{p_M+v}^{-\epsilon_M}\ldots x_{p_1+v}^{-\epsilon_1} = w,
$$
whence the result. \qed

%\begin{Theorem}\label{asphTransfer} Let $\mathcal{Q} = (a,x:a^n,W(a,x))$ be a presentation for a group $E$ where
%
%$$
%W(a,x) = x^{\epsilon_1}a^{p_1}\ldots x^{\epsilon_L}a^{p_L}
%$$
%is such that $L \geq 1$, $0 \leq p_i < n$, and $\epsilon_i = \pm 1$. Suppose that $\nu^f:E \ra C_n$ is a retraction given by $\nu^f(a) = a$ and $\nu^f(x) = a^f$. Let $M$ be the cellular model of the relative presentation $\mathcal{R} = (C_n,x:W)$ associated to $\mathcal{Q}$. Then, $\mathcal{R}$ is orientable and $\pi_2M = 0$ if and only if $\mathcal{P}_n(\rho^f(W(a,x))$ is orientable and combinatorially aspherical.
%\end{Theorem}

\vspace{.25in}

N.~D.~Gilbert and J.~Howie established a connection between asphericity of relative presentations (in the sense of \cite{BP}) and topological asphericity of two-dimensional cellular models of cyclic presentations in \cite[Lemma 3.1]{GH}. The same argument was adapted to other classes of presentations in \cite{BV,CRS,Spag}. Edjvet and Williams established a complete characterization of topological asphericity for two-dimensional cellular models of the cyclic presentations $\mathcal{P}_n(k,l)$ in \cite[Corollary 4.3]{EdjvetWilliams}, correcting an error in \cite[Theorem 4.3]{CRS}. The next theorem sharpens and broadens the scope of these results in several ways, first by eliminating any restriction on the length or shape of the cyclic relators, second by extending to consideration of the more general concept of combinatorial asphericity for cyclic presentations, and third by establishing its logical equivalence with the homotopy condition $\pi_2M = 0$ for relative presentations. The group-theoretic fallout of the homotopy condition $\pi_2M = 0$ is identical to that which derives from the other more specialized forms of asphericity. See Theorems \ref{OrAsphMain} and \ref{freeAction}.

\begin{Theorem}\label{asphTransfer} Let $M$ be the cellular model of an orientable relative presentation $(C_n,x:W)$ for a group $E$ where $C_n$ is the cyclic group of order $n$ generated by $a$ and
$$
W = x^{\epsilon_1}a^{p_1}\ldots x^{\epsilon_L}a^{p_L} \in C_n \ast \<x\>
$$
is cyclically reduced. Suppose that $\nu^f:E \ra C_n$ is a retraction given by $\nu^f(a) = a$ and $\nu^f(x) = a^f$. Then, $\rho^f(W)$ is cyclically reduced and the cyclic presentation $\mathcal{P}_n(\rho^f(W))$ is orientable. Moreover, $\mathcal{P}_n(\rho^f(W))$ is combinatorially aspherical if and only if $\pi_2M = 0$.
\end{Theorem}

\proof Let $w = \rho^f(W) = x_{u(1)}^{\epsilon_1}\ldots x_{u(L)}^{\epsilon_L}$ where each $u(i)$ is defined modulo $n$. If $w$ is not cyclically reduced, then for some $i$ (defined modulo $L$), we have $\epsilon_{i+1} = -\epsilon_i$ and $u(i+1) \equiv u(i)$. Working modulo $n$, this leads to $v(i) + \epsilon_if + p_i \equiv v(i) + \epsilon_if$
%
%\begin{eqnarray*}
%% \nonumber to remove numbering (before each equation)
%  v(i) + \epsilon_if + p_i &=& v(i+1)  \\
%   &=& u(i+1) - \left(\frac{\epsilon_{i+1}-1}{2}\right) \\
%   &=& u(i) - \left(\frac{-\epsilon_{i}-1}{2}\right) \\
%   &=& v(i) + \left(\frac{\epsilon_{i}-1}{2}\right) + \left(\frac{\epsilon_{i}+1}{2}\right) \\
%   &=& v(i) + \epsilon_if,
%\end{eqnarray*}
so $p_i \equiv 0$ modulo $n$. So if $w$ is not cyclically reduced, then $W \in C_n \ast \<x\>$ is not cyclically reduced.

Now suppose that $\mathcal{P}_n(w)$ is not orientable. By Lemma \ref{orientableShape}, $n = 2m$ is even and $w = u\theta^m(u)^{-1}$ for some linearly reduced word $u$. It follows that $L = 2M$ is even and

$$
x_{u(1)}^{\epsilon_1}\ldots x_{u(M)}^{\epsilon_{M}}x_{u(M+1)}^{\epsilon_{M+1}}\ldots x_{u(2M)}^{\epsilon_{2M}} = x_{u(1)}^{\epsilon_1}\ldots x_{u(M)}^{\epsilon_M}x_{u(M)+m}^{-\epsilon_M}\ldots x_{u(M)+1}^{-\epsilon_1},
$$
whence
%
%
%\begin{eqnarray*}
%% \nonumber to remove numbering (before each equation)
%  \rho^f(W) &=& x_{u(1)}^{\epsilon_1}\ldots x_{u(M)}^{\epsilon_{M}}x_{u(M+1)}^{\epsilon_{M+1}}\ldots x_{u(2M)}^{\epsilon_{2M}}\\
%  &=& x_{u(1)}^{\epsilon_1}\ldots x_{u(M)}^{\epsilon_M}\left(x_{u(1)+m}^{\epsilon_1}\ldots x_{u(M)+m}^{\epsilon_M}\right)^{-1} \\
%   &=& x_{u(1)}^{\epsilon_1}\ldots x_{u(M)}^{\epsilon_M}x_{u(M)+m}^{-\epsilon_M}\ldots x_{u(M)+1}^{-\epsilon_1}.
%\end{eqnarray*}
%Thus the following conditions are satisfied for $j = 0,\ldots,M-1$:

\begin{eqnarray}
% \nonumber to remove numbering (before each equation)
  u(M+j+1) &\equiv & u(M-j)+m\ \ \mbox{and} \label{eqn:u}\\
  \epsilon_{M+j+1} &=& -\epsilon_{M-j}\label{eqn:e}
\end{eqnarray}
for $j=0,\ldots,M-1$. The goal is to show that $(C_n,x:W)$ is not orientable by showing that $W = Xa^{p_M}X^{-1}a^{p_{2M}}$ in $C_n \ast\<x\>$ where $\left(a^{p_M}\right)^2 = \left(a^{p_{2M}}\right)^2 = 1$ in $C_n$ and $X = x^{\epsilon_1}a^{p_1}\ldots a^{p_{M-1}}x^{\epsilon_M}$. Using equations (\ref{eqn:u}) and (\ref{eqn:e}), one checks that

\begin{equation}\label{eqn:v}
  v(M+j+1) + \epsilon_{M+j+1}f \equiv v(M-j) + m.
\end{equation}
for $j = 0,\dots ,M-1$. Equations (\ref{eqn:v}) and (\ref{eqn:e}) with $j=0$ imply

$$
v(M) + m \equiv v(M+1) + \epsilon_{M+1}f \equiv v(M) + \epsilon_Mf+p_M - \epsilon_{M}f
$$
so $p_M \equiv m$. Thus $\left(a^{p_M}\right)^2 = 1$ in the coefficient group $C_n = C_{2m}$. Equations (\ref{eqn:v}) and (\ref{eqn:e}) with $j=M-1$ lead to

$$
m \equiv v(1)+m \equiv v(2M) + \epsilon_{2M}f \equiv \left(\sum_{i=1}^{2M-1}\epsilon_if+p_i\right)+\epsilon_{2M}f.
$$
By Lemma \ref{Lemma:splitCrit}, $\sum_{i=1}^{2M}\epsilon_if+p_i \equiv 0$ modulo $n$. Thus $p_{2M} \equiv -m \equiv m$ modulo $n$ and so $\left(a^{p_{2M}}\right)^2 = 1$. Now it suffices to show that $p_{M+j} +p_{M-j} \equiv 0$ modulo $n$ for $j=1,\ldots,M-1$. For this, equations (\ref{eqn:v}) and (\ref{eqn:e}) for the indices $j-1$ and $j$ combine as follows:

\begin{eqnarray*}
% \nonumber to remove numbering (before each equation)
  v(M-j) + m &\equiv &  v(M+j+1)+\epsilon_{M+j+1}f\\
   &\equiv & v(M+j) + \epsilon_{M+j}f + p_{M+j} + \epsilon_{M+j+1}f \\
   &\equiv & v(M-j+1) + m + p_{M+j} + \epsilon_{M+j+1}f \\
   &\equiv & v(M-j) + \epsilon_{M-j}f + p_{M-j} + m + p_{M+j} - \epsilon_{M-j}f.
\end{eqnarray*}
This implies that $p_{M+j} +p_{M-j} \equiv 0$ modulo $n$, as desired.

As in the proof of Theorem \ref{splitCrit}, let $p:\widehat{M} \ra M$ be the regular covering projection with $p_\sharp(\pi_1\widehat{M}) = \ker \nu^f$ and let $K = \widehat{M}/\widetilde{D}_n$. The two-skeleton $K^{(2)}$ is the cellular model of the cyclic presentation $\mathcal{P}_n(w)$ and $\pi_2M \cong \pi_2K$. The two-skeleton $M^{(2)}$ is the cellular model of a lift of the relative presentation $(C_n,x:W)$. This lift has the form $\mathcal{Q} = (a,x:a^n,W(a,x))$ and can be chosen carefully so that $W(a,x)$ has exponent $e$ in the free group with basis $\{a,x\}$. To complete the proof of the theorem requires an understanding of how length two identity sequences for $\mathcal{P}_n(w) = \mathcal{P}_n(\rho^f(W))$ are related to three-cells of $K$.

If $W$ is not a proper power in $C_n \ast \<x\>$ then $K = K^{(2)}$ is two-dimensional. If, on the other hand, $W = \dot{W}^e$ has exponent $e > 1$ in $C_n \ast \<x\>$, then the three-cell of the complement $M \backslash D_n$ is attached by a spherical map $S^2 \ra M^{(2)}$ associated to an identity sequence $(1,W(a,x))(\dot{W}(a,x),W(a,x))^{-1}$ for $\mathcal{Q}$. Let $\nu^f(\dot{W}) = a^\sigma \in C_n$ be the image of the root $\dot{W}$ under the retraction $\nu^f$. Further, let $\hat{w} = \hat{w}(x_0,\ldots, x_{n-1}) = \rho^f(\dot{W})$ be the rewrite of the root of $W$. It is worth noting that if $\sigma \not \equiv 0$ modulo $n$ (that is, if the root $\dot{W}$ is not in the kernel of the retraction $\nu^f$), then the word $\hat{w}$ is \emph{not} the root of $w$. Indeed, the exponent of $w = \rho^f(W)$ is equal to the order of the subgroup that is generated by the image $\nu^f(\dot{W}) = a^\sigma \in C_n$ in the cyclic group of order $n$.

The attaching map for the three-cell of $M$ admits $n$ distinct lifts through the covering projection $\widehat{M} \ra M$. These descend to the quotient complex $K$ to provide attaching maps for the three-cells of $K$. These attaching maps correspond to the shifts of the identity sequence

$$
\Delta^f = (1,w)(\hat{w},\theta^\sigma(w))^{-1}
$$
for $\mathcal{P}_n(w)$, that is, to

%$$
%\Delta^f_k = (1,\theta^k(w))(\theta^k(\hat{w}),\theta^{k+\sigma}(w))^{-1}
%$$
$$
\theta^k\Delta^f = (1,\theta^k(w))(\theta^k(\hat{w}),\theta^{k+\sigma}(w))^{-1}
$$
where $k = 0,\ldots, n-1$. Each of these identity sequences determines a homotopy class in the kernel of the inclusion-induced homomorphism $\pi_2K^{(2)} \ra \pi_2K$.

If $\pi_2M = 0$, then $\pi_2K = 0$, which implies that $\pi_2K^{(2)}$ is generated by the classes of these length two identity sequences for $\mathcal{P}_n(w)$, whence $\mathcal{P}_n(w)$ is combinatorially aspherical by Proposition \ref{CA}.

For the converse, it is necessary and sufficient to show that for each length two identity sequence of the presentation $\mathcal{P}_n(w)$, the corresponding homotopy class is trivial in $\pi_2K$. A length two identity sequence for $\mathcal{P}_n(w)$ can arise only if a cyclic permutation of $w$ is identically equal to one of its shifts $\theta^v(w)$. It suffices to consider an identity sequence of the form $\Sigma = (1,w)(s,\theta^v)^{-1}$ where $s$ is an initial segment of $w$. This means that $w = st$ where $s = x_{u(1)}^{\epsilon_1} \ldots x_{u(q)}^{\epsilon_q}$ and $t = x_{u(q+1)}^{\epsilon_{q+1}} \ldots x_{u(L)}^{\epsilon_L}$, and that

$$
ts = x_{u(q+1)}^{\epsilon_{q+1}} \ldots x_{u(L)}^{\epsilon_{L}}x_{u(1)}^{\epsilon_{1}} \ldots x_{u(q)}^{\epsilon_{q}} = x_{u(1)+v}^{\epsilon_{1}} \ldots x_{u(L)+v}^{\epsilon_{L}}
$$
involving cyclically reduced words in the free group with basis $\{x_0,\ldots x_{n-1}\}$. This relation leads to the conditions

\begin{eqnarray}
% \nonumber to remove numbering (before each equation)
  u(q+i) & \equiv & u(i)+v\ \ \mbox{modulo $n$ and} \label{eqn:u2}\\
  \epsilon_{q+i} &=& \epsilon_{i}\label{eqn:e2}
\end{eqnarray}
for $i = 1,\ldots,L$ modulo $L$. In turn, these lead to

\begin{equation}
  v(q+i) \equiv v(i)+v\ \ \mbox{modulo $n$}\label{eqn:v2}
\end{equation}
and

\begin{equation}
  p_{q+i} \equiv p_i\ \ \mbox{modulo $n$}\label{eqn:p2}
\end{equation}
for $i = 1,\ldots,L$. The conditions (\ref{eqn:e2}) and (\ref{eqn:p2}) imply the exponent $e$ of $W$ is at least $e \geq L/\gcd(L,q)$, so the result holds if $W$ is not a proper power in $C_n \ast \<x\>$. Further, if $S = x^{\epsilon_1}a^{p_1}\ldots x^{\epsilon_q}a^{p_q}$ is the initial segment of $W \in C_n \ast \<x\>$ corresponding to $s$, so that $s = \rho^f(S)$, then $S = \dot{W}^\ell$ for suitable $\ell \geq 1$ depending on $e$, $L$, and $q$. (In other words, $S$ commutes with $W$ in $C_n \ast \<x\>$ and so $S$ is a power of the root $\dot{W}$.) Finally, the condition $(\ref{eqn:v2})$ with $i = 1$ implies that $v(q+1) \equiv \nu^f(S) \equiv \ell \sigma \equiv v$ modulo $n$.

Form the following product of identity sequences for $\mathcal{P}_n(w) = \mathcal{P}_n(\rho^f(W))$,

%$$
%\Pi = \prod_{j=0}^{\ell - 1} \left(\left(\prod_{k=0}^{j} \theta^{k\sigma}(\hat{w})\right)\Delta^f_{j\sigma}\right).
%$$

$$
\Pi = \prod_{j=0}^{\ell - 1} \left(\left(\prod_{k=0}^{j-1} \theta^{k\sigma}(\hat{w})\right) \cdot \theta^{j\sigma}\Delta^f\right),
$$
and note that $[\Pi ] \in \ker (\pi_2K^{(2)} \ra \pi_2K)$. Expanding terms, the product $\Pi$ telescopes as follows:

\begin{eqnarray*}
\Pi &=& \left(\Delta^f\right) \left(\hat{w}\cdot \theta^\sigma \Delta^f\right) \left(\hat{w}\theta^\sigma(\hat{w}) \cdot \theta^{2\sigma}\Delta^f\right) \cdots \left(\hat{w}\theta^\sigma(\hat{w})\cdots \theta^{\sigma(\ell - 2)}(\hat{w})\cdot \theta^{\sigma(\ell -1)}\Delta^f\right)\\
&=& \left(1,w\right)\left(\hat{w},\theta^\sigma(w)\right)^{-1} \cdot \left(\hat{w},\theta^\sigma(w)\right)\left(\hat{w}\theta^\sigma(\hat{w}),\theta^{2\sigma}(w)\right)^{-1}\\
&& \cdot \left(\hat{w}\theta^\sigma(\hat{w}),\theta^{2\sigma}(w)\right)\left(\hat{w}\theta^\sigma(\hat{w})\theta^{2\sigma}(\hat{w}),\theta^{3\sigma}(w)\right)^{-1} \\
&& \cdots \\
&& \cdot \left(\hat{w}\theta^\sigma(\hat{w})\cdots \theta^{\sigma(\ell - 2)}(\hat{w}),\theta^{\sigma(\ell -1)}(w)\right)\left(\hat{w}\theta^\sigma(\hat{w})\cdots \theta^{\sigma(\ell - 1)}(\hat{w}),\theta^{\sigma\ell}(w)\right)^{-1}\\
&=& (1,w)\left(\hat{w}\theta^\sigma(\hat{w})\cdots \theta^{\sigma(\ell - 1)}(\hat{w}),\theta^{\sigma\ell}(w)\right)^{-1}\\
&=& (1,w) \left(\hat{w}\theta^\sigma(\hat{w})\cdots \theta^{\sigma(\ell - 1)}(\hat{w}),\theta^v(w)\right)^{-1}.
\end{eqnarray*}
The argument is concluded by noting that since $\rho^f(\dot{W}) = \hat{w}$ and $\nu^f(\dot{W}) = a^\sigma$, whence
\begin{eqnarray*}
  \hat{w}\theta^\sigma(\hat{w})\cdots \theta^{\sigma(\ell - 1)}(\hat{w}) &=& \rho^f(\dot{W})\theta^\sigma(\rho^f(\dot{W}))\cdots \theta^{\sigma(\ell - 1)}(\rho^f(\dot{W})) \\
  &=& \rho^f(\dot{W}^\ell) = \rho^f(S) = s
\end{eqnarray*}
and so $\Pi = \Sigma$, whence $[\Sigma] \in \ker (\pi_2K^{(2)} \ra \pi_2K)$. \qed

\vspace{.25in}

To prove Theorem A, consider an orientable and combinatorially aspherical cyclic presentation $\mathcal{P}_n(w)$ where $w = w(x_0,\ldots ,x_{n-1})$ is nonempty and cyclically reduced. Write

$$
w = x_{p_1}^{\epsilon_1}\ldots x_{p_L}^{\epsilon_L}
$$
where $L \geq 1$, $0 \leq p_j < n$, and $\epsilon_j = \pm 1$. After rewriting $x_j = a^jxa^{-j}$ and cyclically permuting, let

$$
W = x^{\epsilon_1}a^{p_2-p_1}x^{\epsilon_2}a^{p_3-p_2} \cdots x^{\epsilon_L}a^{p_1-p_L}
$$
and consider the relative presentation $(C_n,x:W)$ for the group $E = (C_n \ast \<x\>)/\<\<W\>\>$. There is a retraction $\nu^0:E \ra C_n$ and Theorem \ref{splitCrit} shows that the kernel of $\nu^0$ is given by the cyclic presentation $\mathcal{P}_n(\rho^0(W)) = \mathcal{P}_n(w)$. (Note that $\theta^{p_1}(\rho^0(W)) = w$.) It is shown below that the relative presentation $(C_n,x:W)$ is orientable. This being the case, Theorem \ref{asphTransfer} implies that the cellular model $M$ of $(C_n,x:W)$ has $\pi_2M = 0$ (and so $M$ is aspherical by Proposition \ref{relComb}(c)). That $C_n$ acts freely via the shift on the nonidentity elements of $G_n(\rho^0(W)) = G_n(w)$ then follows from Theorem \ref{freeAction}.

For orientability, suppose that $(C_n,x:W)$ is not orientable. This implies that $W = X\alpha X^{-1}\beta \in C_n \ast \<x\>$ where $\alpha^2 = \beta^2 = 1$ in $C_n$. Set $q_j = p_{j+1}-p_j$ for $j = 1,\ldots, L$ modulo $L$. With $X = x^{\epsilon_1}a^{q_1} \cdots a^{q_{M-1}}x^{\epsilon_M}$, it follows that $L = 2M$, and the following conditions are satisfied:

\begin{eqnarray*}
% \nonumber to remove numbering (before each equation)
 \left(a^{q_M}\right)^2 &=& \left(a^{q_{2M}}\right)^2 = 1, \\
 \epsilon_{M+j} &=& -\epsilon_{M-j+1}\ \ (j = 1,\ldots, M), \\
  q_{M+j} &\equiv & -q_{M-j}\ \ (j = 1,\ldots, M-1).
\end{eqnarray*}
The fact that $w$ is cyclically reduced implies that $n = 2m$ is even and $q_M \equiv q_{2M} \equiv m$ modulo $n$, whence $p_{M+1} \equiv p_M + m$ and $p_{2M} \equiv p_1 + m$. Inductively, the conditions on the $q$'s imply that $p_{M+j} \equiv p_{M-j+1} + m$ for $j = 1,\dots, M$. Letting $u = x_{p_1}^{\epsilon_1}\ldots x_{p_M}^{\epsilon_M}$, this implies that $u\theta^m(u)^{-1} = w$, contrary to the fact that $\mathcal{P}_n(w)$ is orientable. This completes the proof of Theorem A.

\section{Proofs of Theorems B and C}

Now consider the groups $G_n(k,l)$ given by cyclic presentations

$$
\mathcal{P}_n(k,l) = (x_0,\ldots,x_{n-1}\ |\ x_ix_{i+k}x_{i+l}\ (i=0,\ldots,n-1)).
$$
Each group $G_n(k,l)$ has a homomorphic image of order three and so is nontrivial. %We first consider the case where $d=\gcd(n,k,l) > 1$.

\begin{Lemma}\label{CommonFactor}
  Let $d = \gcd(n,k,l)$ and let $\theta$ be the shift on $G_n(k,l)$.
  \begin{enumerate}
    \item[(a)] $\mathcal{P}_n(k,l)$ is a disjoint union of $d$ subpresentations $\mathcal{Q}_1,\ldots, \mathcal{Q}_d$, where $\mathcal{Q}_j$ consists of those generators and relators of $\mathcal{P}_n(k,l)$ that involve only the generators $x_i$ where $i \equiv j \mod d$ and each $Q_j$ is isomorphic to the cyclic presentation $\mathcal{P}_{n/d}(k/d,l/d)$.
    \item[(b)] The shift $\theta$ on $G_n(k,l)$ carries the free factor determined by the subpresentation $\mathcal{Q}_j$ onto the free factor determined by $\mathcal{Q}_{j+1}$ (subscripts modulo $d$) and the shift $\theta'$ on $G_{n/d}(k/d,l/d)$ is the restriction of $\theta^d$.
    \item[(c)] A power $\theta^p$ of the shift on $G_n(k,l)$ has a non-identity fixed point if and only if $d\,|\,p$ and $\theta'^{p/d}$ has a non-identity fixed point on $G_{n/d}(k/d,l/d)$.
%
%    \item[(b)] The shift on $G_n(k,l)$ has an orbit of size $p > 1$ if and only if the shift on $G_{n/d}(k/d,l/d)$ has an orbit of size $p/d$ that does not contain the identity.
  \end{enumerate}
\end{Lemma}

\proof The statement (a) was noted in \cite{CRS,EdjvetWilliams} and (b) follows immediately. To prove (c), if $p\,|\,d$ and $\theta'^{p/d}$ has a non-identity fixed point, then so does $\theta^p$, by (b). Conversely, suppose that $\theta^p$ fixes a nonempty reduced word $1\neq w = a_{j(1)}\ldots a_{j(\ell)}$ in the free product, where each $a_{j(t)}$ is a nontrivial element in one of the free factors of $\mathcal{P}_n(k,l)$ determined by $\mathcal{Q}_{j(t)}$ ($1 \leq j(t) \leq d$). Thus consecutive letters $a_{j(t)}\neq 1$ lie in distinct factors. Now (b) implies that $p\,|\,d$ and $\theta^p(a_{j(1)}) = \theta'^{p/d}(a_{j(1)}) = a_{j(1)}$. \qed

\vspace{.25in}

Using \cite[Theorem 4.2]{CCH}, Lemma \ref{CommonFactor}(a) implies that $\mathcal{P}_n(k,l)$ is combinatorially aspherical if and only if $\mathcal{P}_{n/d}(k/d,l/d)$ is combinatorially aspherical. Lemma \ref{CommonFactor}(c) implies that the shift $\theta$ acts freely on the nonidentity elements of $G_n(k,l)$ if and only if the nonidentity elements of $G_{n/d}(k/d,l/d)$ are permuted freely by its shift $\theta'$. In addition, if $\gcd(n,k,l) > 1$, then Lemma \ref{CommonFactor}(a) implies that $G_n(k,l)$ is infinite and Lemma \ref{CommonFactor}(b) implies that $\theta$ has no nontrivial fixed points. These remarks show that it suffices to prove Theorems B and C in the case where $\gcd(n,k,l) = 1$.

%\vspace{.25in}

Edjvet and Williams developed a taxonomy for the parameters $(n,k,l)$ and used it to characterize asphericity for the cellular models of the presentations $\mathcal{P}_n(k,l)$ \cite[Theorem A]{EdjvetWilliams} and finiteness for the groups $G_n(k,l)$ \cite[Theorem B]{EdjvetWilliams}. Their taxonomy involved the following divisibility conditions.

\begin{enumerate}
  \item[(A)] $3\mid n$ and $3\mid k+l$.
  \item[(B)] $n\mid k+l$ or $n\mid 2k-l$ or $n\mid 2l-k$.
  \item[(C)] $n\mid 3k$ or $n\mid 3l$ or $n\mid 3(k-l)$.
\end{enumerate}

\noindent Assuming that $\gcd(n,k,l) = 1$, the strategy is to show that Theorems B and C hold when condition (B) is satisfied, and then again when condition (C) is satisfied. The case when conditions (B) and (C) do not hold is then treated separately. The semidirect product $E_n(k,l) = G_n(k,l) \rtimes_\theta C_n$ has presentation $\mathcal{Q}_n(k,l) = (a,x:a^n, W(a,x))$ where

$$
W(a,x) = xa^kxa^{l-k}xa^{-l}.
$$
Viewing $W \in C_n \ast \<x\>$, the relative presentation $(C_n,x:W)$ is orientable.

\begin{Lemma}\label{lemma:B}
 Assume that $\gcd(n,k,l) = 1$ and that the condition (B) is satisfied. With a change of variable of the form $u = xa^t$ for suitable $t$, the relator $W$ is conjugate in $C_n \ast \<x\> = C_n \ast \<u\>$ to a word of the form $u^3a^{3p}$, where $\gcd(p,n) = 1$. Consequently:
  \begin{enumerate}
    \item[(a)] The element $a^3$ is central in the semidirect product $E_n(k,l)$ and so the shift $\theta$ on $G_n(k,l)$ satisfies $\theta^3 = 1$.
    \item[(b)] The cyclic presentation $\mathcal{P}_n(k,l)$ is combinatorially aspherical if and only if $n=3$.
    \item[(c)] If $n =3$, then $E_3(k,l) \cong C_3 \ast C_3$ and $C_3$ acts freely via the shift on the nonidentity elements of $G_3(k,l)$.
    \item[(d)] If $3\mid n$ and $n\neq 3$, then $E_n(k,l) \cong \<u\> \ast_{C_{n/3}} C_n$ is infinite and the shift $\theta$ has no fixed points.
    \item[(e)] If $3\nmid n$ then $E_n(k,l) \cong C_{3m}$ is finite cyclic, so $\theta = 1$.
  \end{enumerate}
\end{Lemma}

\proof Taking cases in condition (B), if $n\mid k+l$, then let $u = xa^k$ and rewrite $W = u^2a^{l-2k}ua^{-(k+l)}$, which cyclically permutes to $ u^3a^{-3k}$ in $C_n \ast\<x\> = C_n \ast \<u\>$. The conditions $\gcd(n,k,l) = 1$ and $n\mid k+l$ imply that $\gcd(n,-k) = 1$. The other cases in condition (B) proceed similarly. The conclusion (a) follows immediately, so if $n \neq 3$, then $C_n$ does not act freely via the shift, whence Theorem A implies that $\mathcal{P}_n(k,l)$ is not combinatorially aspherical. If $n=3$, then the relative presentation $(C_n,u:u^3a^{3p}) = (C_n,u:u^3)$ is aspherical in the sense of \cite{BP}, so $\mathcal{P}_3(k,l)$ is combinatorially aspherical by Theorem \ref{asphTransfer}. If $3\mid n$ and $n\neq 3$ then $E \cong \<u\> \ast_{u^3=a^{-3p}} C_n$ splits as a nontrivial free product with amalgamation where the element $a \in C_n$ lies outside the amalgamating subgroup of order $n/3$ generated by $a^3$. Thus $a$ generates its centralizer in $E$ \cite[Theorem 4.5]{MKS} and so no nonidentity element of $G_n(k,l) \leq E$ is centralized by $a$. The remaining conclusions are clear. \qed

%\vspace{.25in}

\paragraph{Condition (B):} The conclusions of Lemma \ref{lemma:B} show that Theorems B and C both hold true when $\gcd(n,k,l) = 1$ and the condition (B) is satisfied. That is:

\begin{itemize}
  \item $G_n(k,l)$ is finite  $\iff 3 \nmid n \iff \theta = 1 \iff \theta$ has a nonidentity fixed point; and
  \item $\mathcal{P}_n(k,l)$ is combinatorially aspherical $\iff n=3 \iff C_n$ acts freely via the shift on the nonidentity elements of $G_n(k,l)$.
\end{itemize}

Edjvet and Williams \cite{EdjvetWilliams} showed that if $3\mid n$ in the situation of Lemma \ref{lemma:B}, then $G_n(k,l)$ is isomorphic to the free group of rank two. Exemplars for this case are $\mathcal{P}_{3k}(1,2) \cong F(x_0,x_1)$ where the shift is given by $\theta(x_0) = x_1$ and $\theta(x_1)= x_2 = (x_0x_1)^{-1}$, under which $C_3$ acts freely on the nonidentity elements of the free group $F(x_0,x_1)$. The next lemma describes a very different set of exemplars that apply when the condition (C) is satisfied.

\begin{Lemma}\label{lemma:cyc}
  If $\gcd(p,n) = 1$, then the cyclic presentation $\mathcal{P}_n(0,p)$ defines a cyclic group $G_n(0,p) \cong C_s$ of order $s = 2^n - (-1)^n$. Moreover the shift $\theta$ fixes a subgroup of order three.
\end{Lemma}

\proof Taking subscripts modulo $n$, the relations combine sequentially to show:

$$
x_i = x_{i-p}^{-2} = x_{i-2p}^4 = \cdots = x_{i-kp}^{(-2)^k}\ \ \mbox{for all $i$ and $k$}.
$$
Since $\gcd(p,n)=1$, this shows that $G_n(0,p)$ is cyclic of order $s = 2^n - (-1)^n$, generated by $x_0$. The shift $\theta$ satisfies $\theta^p(x_0) = x_p = x_0^{-2}$ and therefore $\theta^p\left(x_0^{s/3}\right) = x_0^{-2s/3} = x_0^{s/3}$. This means that $\theta^p$, and hence $\theta$, fixes a subgroup of order three. \qed

\begin{Lemma}\label{lemma:C}
  Assume that $\gcd(n,k,l) = 1$ and that the condition (C) is satisfied.
   \begin{enumerate}
     \item[(a)] If $3 \nmid n$, then $\mathcal{P}_n(k,l) = \mathcal{P}_n(0,p)$ for some $p$ satisfying $\gcd(p,n) = 1$.
     \item[(b)] If $3\mid n$ and $3 \nmid k+l$, then for some $f$ and for some $p$ satisfying $\gcd(p,n) = 1$ and $n \mid 3f$, the semidirect product $E = E_n(k,l)$ admits a retraction $\nu^f: E \ra C_n$ where $\mathcal{P}_n(\rho^f(W)) = \mathcal{P}_n(0,p)$.
     \item[(c)] If the condition (A) is satisfied, so that $3 \mid n$ and $3 \mid k+l$, then for some $f$ and for some $p$ satisfying $\gcd(p,n) = 3$ and $n \mid 3f$, the semidirect product $E = E_n(k,l)$ admits a retraction $\nu^f: E \ra C_n$ where $\mathcal{P}_n(\rho^f(W)) = \mathcal{P}_n(0,p)$.
   \end{enumerate}
\end{Lemma}

\proof Suppose first that $3 \nmid n$. Taking cases in condition (C), if $n \mid 3k$, then $n \mid k$, so $\mathcal{P}_n(k,l) = \mathcal{P}_n(0,l)$ and $\gcd(n,l) = 1$. The other cases are handled in the same way, using the fact that $\mathcal{P}_n(k,0) = \mathcal{P}_n(0,k)$ and $\mathcal{P}_n(k,k) = \mathcal{P}_n(0,n-k)$.

Now suppose that $n$ is divisible by three. Again taking cases in condition (C), if $n\mid 3k$, then we have the retraction $\nu^{-k}: E \ra C_n$ where $\rho^{-k}(W) = x_0^2x_{l-2k}$. Using the fact that $\gcd(n,k,l) = 1$, if $3 \nmid k+l$, then $\gcd(l-2k,n) = 1$, whereas if $3 \mid k+l$, then $\gcd(l-2k,n) = 3$. The other cases in condition (C) are handled similarly; if $n\mid 3l$, use the retraction $\nu^l$, and if $n\mid 3(k-l)$, then use the retraction $\nu^{k-l}$. \qed

%\vspace{.25in}

\paragraph{Condition (C):} Lemmas \ref{TransferDynamic}, \ref{CommonFactor}, \ref{lemma:cyc}, and \ref{lemma:C} combine to show that Theorem C holds true when $\gcd(n,k,l) = 1$ and the condition (C) is satisfied. That is:

\begin{itemize}
  \item $G_n(k,l)$ is finite  $\iff$ condition (A) is not satisfied $\iff \theta$ has a nonidentity fixed point.
\end{itemize}
Theorem B also holds here because the shift determines a free action on the nonidentity elements of $G_n(k,l)$ only when $n=3 \mid k+l$, in which case Lemma \ref{lemma:B} applies. In all other cases, the cyclic presentation $\mathcal{P}_n(k,l)$ is not combinatorially aspherical by Theorems \ref{OrAsphMain}(c) and \ref{asphTransfer} because $G_n(0,p)$ has nontrivial torsion.

\vspace{.25in}

Lemma \ref{lemma:C} provides quick access to structural results that were obtained through detailed analysis in \cite[Lemma 2.5, Lemma 3.2]{EdjvetWilliams}. For example, if $(n,k,l)$ satisfy the hypotheses of Lemma \ref{lemma:C}, then $G_n(k,l)$ is virtually a free product of finite cyclic groups. As the next lemma illustrates, commensurability provides additional information.

\begin{Lemma}\emph{\textbf{(\cite[Lemma 3.2]{EdjvetWilliams})}}\label{lemma:metacyclic} Assume that $\gcd(n,k,l) = 1$ and that condition (C) is satisfied. If $3\mid n$ and $3 \nmid k+l$ (so that condition (A) is not satisfied) then $G_n(k,l)$ is metacyclic of order $s = 2^n - (-1)^n$.
\end{Lemma}

\proof By Lemma \ref{lemma:C}(b), there is a retraction $\nu^f: E \ra C_n$ where $n \mid 3f$ and $\ker \nu^f \cong G_n(0,p)$ is such that $\gcd(n,p) = 1$. By Lemma \ref{lemma:cyc}, $\ker \nu^f$ is cyclic of order $s$, and is generated by $xa^{-f}$ by Theorem \ref{splitCrit}. Now $G_n(k,l) = \ker \nu^0$ and $\nu^0((xa^{-f})^3) = 1$ so $\ker \nu^0 \cap \ker \nu^f$ has index three in $C_s = \ker \nu^f$ and the result follows. \qed

\vspace{.25in}

To conclude the proofs of Theorems B and C, it remains to consider the cases where neither (B) nor (C) is satisfied, at which point asphericity takes over.

\begin{Theorem}\emph{\textbf{(\cite[Theorem A]{EdjvetWilliams})}}
  Assume that $\gcd(n,k,l) = 1$. If neither the condition (B) nor the condition (C) is satisfied, then the two-dimensional cellular model of $\mathcal{P}_n(k,l)$ is aspherical (in the topological sense), unless $n = 18$ and $3 \mid k+l$.
\end{Theorem}

\noindent As was noted earlier, if the cellular model of $\mathcal{P}_n(k,l)$ is aspherical in the topological sense, then the presentation is combinatorially aspherical and so the shift acts freely on the nonidentity elements by our Theorem A. Further, in this case the group $G_n(k,l)$ is torsion-free and so, being nontrivial, is infinite. Therefore it remains to show that Theorems B and C are true in the case where $n = 18$ and $3 \mid k+l$. In this case, Edjvet and Williams showed that $G_{18}(k,l)$ is the free product of a free group of rank two with the cyclic group of order $19$ \cite[page 771]{EdjvetWilliams}. Thus $G_{18}(k,l)$ is infinite and $\mathcal{P}_{18}(k,l)$ is not combinatorially aspherical. It remains to show that the shift $\theta$ has no nonidentity fixed point in $G_{18}(k,l)$, but that some nonidentity power of $\theta$ has a nonidentity fixed point in $G_{18}(k,l)$. %For this, we need to investigate $E_{18}(k,l) = G_{18}(k,l) \rtimes_\theta C_{18}$.

\begin{Lemma}\label{lemma:18}
  Assume that $n=18$, that $\gcd(18,k,l) = 1$, that $3 \mid k+l$, and that neither condition (B) nor condition (C) is satisfied. With $E_{18}(k,l) = G_{18}(k,l) \rtimes_\theta C_{18}$ defined by the presentation $(a,x:a^{18},xa^kx^{l-k}xa^{-l})$, the split short exact sequence

  $$
  1 \ra G_{18}(k,l) \ra E_{18}(k,l) \stackrel{\nu^0}{\ra} C_{18} \ra 1
  $$
  is equivalent to one of the form

  $$
  1 \ra G \ra E \stackrel{\nu}{\ra} C_{18} \ra 1
  $$
  where $E$ has presentation $(a,u:a^{18},u^2a^9ua^6)$ and $\nu:E \ra C_{18}$ is a retraction onto the cyclic group of order $18$ generated by $a$.
\end{Lemma}

\proof Applying the change of variable $u = xa^{-l}$ to the presentation $(a,x:a^{18},xa^kx^{l-k}xa^{-l})$ leads to $(a,u:a^{18},u^2a^{k+l}u^{2l-k})$. Both $k+l$ and $2l-k$ are defined modulo $18$ and are divisible by $3$. Working modulo $18$, $k+l \not \equiv \pm(2l-k)$ and $k+l,2l-k \not \equiv 0$ because neither condition (B) nor (C) is satisfied. And because $\gcd(n,k,l) = 1$, it follows that $(k+l)+(2l-k) \not \equiv 9$. All of this means that

$$
\{k+l,2l-k\} \equiv \{3,9\}, \{3, 12\}, \{6,9\}, \{6,15\}, \{9,12\},\ \mbox{or}\ \{9,15\}.
$$
%Inverting the generators $a$ and $u$ of $\mathcal{Q}$ leads to $(a,u:a^{-18},(u^2a^{2l-k}ua^{k+l})^{-1})$ so the parameters $k+l$ and $2l-k$ can be taken in either order without affecting the extension class.
All of these choices lead to the same extension. To see this, begin with the parameters $(k,l) = (1,11)$, which satisfy the hypotheses of the lemma and determine the semidirect product $E_{18}(1,11)$ with presentation $(a,x:a^{18},xaxa^{10}xa^{-11})$. The following changes of variable do not affect the group or extension class:
\begin{eqnarray*}
% \nonumber to remove numbering (before each equation)
  u=xa &\ra& (a,u:a^{18},u^2a^9ua^{-12}) \cong (a,u:a^{18},u^2a^9ua^6) \\
  u=xa^{10} &\ra& (a,u:a^{18},u^2a^{-21}ua^{-9}) \cong (a,u:a^{18},u^2a^{15}ua^9) \\
  u=xa^{-11} &\ra& (a,u:a^{18},u^2a^{12}ua^{21}) \cong (a,u:a^{18},u^2a^{12}ua^3)
\end{eqnarray*}
Replacing a relator $u^2a^pua^q$ with its inverse and cyclically permuting in $C_{18} \ast \<u\>$ leads to $\left(u^2a^{18-q}ua^{18-p}\right)^{-1}$. Replacing $a$ and $u$ by their inverses in $u^2a^pua^q$ and cyclically permuting in $C_{18} \ast \<u\>$ leads to $\left(u^2a^qua^p\right)^{-1}$. With $(k,l)$ given as above and letting $E \cong E_{18}(1,11)$ be the group with presentation $(a,u:a^{18},u^2a^9ua^6)$, this implies that there is an isomorphism $\psi: E_{18}(k,l) \ra E$ that carries $a$ to $a^\epsilon$ for $\epsilon = \pm 1$. Now the retraction $\nu:E \ra C_{18}$ given by the composite $E \stackrel{\psi^{-1}}{\ra} E_{18}(k,l) \stackrel{\nu^0}{\ra} C_{18} \stackrel{\epsilon}{\ra} C_{18}$ satisfies $\nu(a) = a$. \qed

\vspace{.25in}

Let $\mathcal{P}_{18}(k,l)$ be a cyclic presentation where $\gcd(18,k,l) = 1$, where $3 \mid k+l$, and where neither (B) nor (C) is satisfied. By Lemma \ref{lemma:18}, the cyclically presented group $G_{18}(k,l)$ can be obtained as the kernel of a retraction $\nu:E \ra C_{18}$ where $E$ has presentation $(a,u:a^{18},u^2a^9ua^6)$, with the shift on $G_{18}(k,l)$ (or its inverse) given by conjugation by $a \in E$. Now $E$ admits retractions $\nu^f:E \ra C_{18}$ where $\nu^f(a) = a$ and $\nu^f(u) = a^f$ for $f = 1,7,13$. These retractions differ by automorphisms of $C_{18}$ and so all have the same kernel, which admits distinct cyclic presentations: $\ker \nu^1 = G_{18}(1,11) = \ker \nu^7 = G_{18}(7,5) = \ker \nu^{13} = G_{18}(1,14) \subseteq E$.

Now consider the group $K$ with presentation $(b,u:b^6, u^2b^3ub^2)$, which appeared in \cite[page 25]{BP}. The element $b$ has order $6$ in $K$ \cite{Levin} and a computer-assisted coset enumeration, for example using GAP \cite{GAP}, shows that $K$ is finite of order $342$. The group $E$ decomposes as a free product with amalgamation

$$
E = K \ast_{b=a^3} C_{18} = K \ast_{C_6} C_{18}.
$$
Since $C_{18} \cap \ker \nu = 1$ and $a$ lies outside the amalgamating subgroup $C_6$ generated by $b = a^3$, the shift on $G_{18}(k,l)$ has no fixed points, as in Theorem C. However, the subgroup $C_6 \leq K$, generated by $b$, acts by left multiplication on the left coset space $K/C_6$. Since $|K/C_6| = 57$, the $C_6$-action on $K/C_6 - \{1C_6\}$ is not free. In fact, a GAP calculation \cite{GAP} shows that the $C_6$-action on $K/C_6$ actually has three fixed points, including $1C_6$ as well as two others. This means that $g \in K \backslash C_6$ with $a^3 g = bg \in gC_6 \subseteq gC_{18}$, whence $g \in E\backslash C_{18}$ and $a^3gC_{18} = gC_{18}$. By Lemma \ref{TransferDynamic}, $\theta^3$ has a nonidentity fixed point on $G_{18}(k,l)$, which concludes the proof of Theorem B.

%$K$ abelianizes to $C_{18}$, it follows that $\ker \nu \cap K = [K,K]$ is cyclic of order $19$. In addition, the element $a \in E$ generates its centralizer in $E$, so

\noindent {\small Department of Mathematics, Kidder 368, Oregon State University, Corvallis, OR 97331 USA\\
bogley@math.oregonstate.edu} %and Nicole Seaders\footnote{About to enter the job market}}

\end{document}